\documentclass[12pt]{amsart}
\usepackage{amsmath,amsthm,amsfonts,amscd,amssymb,amsopn,color}
\usepackage{eucal, mathrsfs}
\usepackage[all]{xy}
\usepackage{fullpage}

\theoremstyle{plain}

\newcommand{\mbb}[1]{\mathbb #1}

\newcommand{\mc}[1]{\mathcal #1}

\newcommand{\ra}{\rightarrow}
\newcommand{\hra}{\hookrightarrow}
\newcommand{\dra}{\dashrightarrow}

\newcommand{\del}{\partial}

\newcommand{\Pic}[1]{\operatorname{Pic}(#1)}
\newcommand{\Picnot}[1]{\operatorname{Pic}^0(#1)}
\newcommand{\sPic}[1]{\operatorname{\mathbf{Pic}}(#1)}
\newcommand{\sPicnot}[1]{\operatorname{\mathbf{Pic}}^0(#1)}
\newcommand{\sPicn}[2]{\operatorname{\mathbf{Pic}}^{#1}(#2)}
\newcommand{\Hgal}{{H}}
\newcommand{\Het}{{H}}

\DeclareMathOperator{\Br}{Br}

\DeclareMathOperator{\ram}{ram}
\DeclareMathOperator{\ind}{ind}

\DeclareMathOperator{\per}{per}

\DeclareMathOperator{\Div}{Div}
\DeclareMathOperator{\NS}{NS}
\DeclareMathOperator{\Prin}{Prin}
\DeclareMathOperator{\im}{im}
\DeclareMathOperator{\coker}{coker}
\DeclareMathOperator{\Spec}{Spec}

\newcommand{\blank}{\ \_ \ }

\newcommand{\sep}[1]{{{#1}^{sep}}}

\newtheorem{thm}{Theorem}[section]
\newtheorem{Lem}[thm]{Lemma}
\newtheorem{cor}[thm]{Corollary}
\newtheorem{Prop}[thm]{Proposition}
\newtheorem{Ex}[thm]{Example}
\newtheorem{defn/Prop}[thm]{Definition/Proposition}
\newtheorem{Remark}[thm]{Remark}

\newtheorem*{thm*}{Theorem}
\newtheorem*{Rem*}{Remark}
\newtheorem*{Lem*}{Lemma}

\newtheorem*{cor*}{Corollary}
\newtheorem*{Prop*}{Proposition}
\theoremstyle{definition}
\newtheorem{defn}[thm]{Definition}

\newenvironment{Rem}{\begin{Remark}}{\nolinebreak{\hfill
$\blacklozenge$} \end{Remark}}

\newcommand{\ov}{\overline}
\newcommand{\til}{\widetilde}

\newcommand{\bZ}{\mathbb{Z}}
\newcommand{\bP}{\mathbb{P}}
\newcommand{\bA}{\mathbb{A}}
\newcommand{\bQ}{\mathbb{Q}}
\newcommand{\bG}{\mathbb{G}}

\newcommand{\E}{{E}}

\usepackage[OT2,T1]{fontenc}
\DeclareSymbolFont{cyrletters}{OT2}{wncyr}{m}{n}
\DeclareMathSymbol{\Sha}{\mathalpha}{cyrletters}{"58}

\newcommand{\fa}{\mathfrak{a}}

\def\<{\left<}
\def\>{\right>}

\begin{document}
\author{Mirela Ciperiani and Daniel Krashen}
\address{Mirela Ciperiani \\
University of Texas at Austin \\
Austin, TX\\
USA}
\email{mirela@math.utexas.edu}
\address{Daniel Krashen \\
University of Georgia, Athens \\
Athens, GA\\
USA}
\email{dkrashen@math.uga.edu}

\title{Relative Brauer groups of genus 1 curves}

\thanks{The authors would like to thank Brian Conrad and Adrian Wadsworth
for helpful conversations and suggestions during the writing of this
paper. Thanks also to Peter Clark for many valuable suggestions and
corrections. Finally we would like to thank the anonymous referee for
a number of helpful comments and corrections.}

\begin{abstract}
In this paper we develop techniques for computing the relative Brauer
group of curves, focusing particularly on the case where the genus is
$1$. We use these techniques to show that the relative Brauer group may
be infinite (for certain ground fields) as well as to determine this
group explicitly for certain curves defined over the rational numbers.
To connect to previous descriptions of relative Brauer groups in the
literature, we describe a family of genus $1$ curves, which we call
``cyclic type'' for which the relative Brauer group can be shown to have
a particularly nice description. In order to do this, we discuss a
number of formulations of the pairing between the points on an elliptic
curve and its Weil-Chat\^elet group into the Brauer group of the ground
field, and draw connections to the period-index problem for genus $1$
curves.
\end{abstract}

\maketitle

%\tableofcontents

\section{Introduction}

Let $X$ be a smooth projective curve of genus $1$ over a field $k$. The
main objects of study in this paper are the elements of the Brauer group
of $k$ which are split by the function field of $X$ i.e., the kernel of
the natural homomorphism $\Br(k) \to \Br(k(X))$. These elements form a
subgroup which we call the relative Brauer group of $X$. This group is of
interest both from the point of view of studying the curve $X$ as well as
from the point of view of field arithmetic and the structure of division
algebras. From the perspective of the curve $X$, one may interpret the
relative Brauer group as an obstruction to the existence of a rational
point, related to the so called ``elementary obstruction.'' In addition,
it is closely related to the period-index problem. From the point of view
of field arithmetic and division algebras, this type of splitting
information for function fields, and more generally index reduction
formulas as in \cite{MPW}, play an important role in constructing examples
and counterexamples (such as Merkurjev's construction of fields with
various $u$-invariants \cite{Lam:Merk-u}).  Unfortunately, such information is
only known for very special varieties at this point, such as for
projective homogeneous varieties under a linear algebraic group. In
particular, in the case of curves, until recently, one only had a complete
description of the relative Brauer group when the genus was $0$. In
\cite{Han:RB}, the relative Brauer groups of certain genus $1$
hyperelliptic curves was described in a surprisingly tractable way. 

In this paper we introduce tools for computing these relative Brauer
groups, and we introduce the notion of a curve of ``cyclic type'' in
order to explain when particularly nice descriptions of the relative
Brauer group may be given as in \cite{Han:RB}. The present paper has
been used in \cite{HHW} to obtain explicit descriptions of the
relative Brauer group for certain plane cubic curves, and in
\cite{Kuo}, which uses generalized Clifford algebra constructions to
study a somewhat more general class of curves.  These tools are also
applied to the period-index problem and the elementary obstruction.

A main ingredient in this paper is a number of reformulations of the
pairing of Tate relating to points on an elliptic curve and the
elements of its Tate-Shafarevich group. Some of these results appear
in the literature, and we cite results of Bashmakov \cite{Ba} and
Lichtenbaum \cite{Lich} to obtain some of these. Although we don't
make the claim that the other reformulations in this paper are
necessarily new, we were unable to obtain references in the literature
after talking to experts in the area, and in any case, we believe that
it is of value to collect some of these results together as we do
here.

The contents of this paper are as follows. Section~\ref{pic period index} presents the basic relationship
between the relative Brauer group, the Picard variety and the period-index problem in terms of a surjective
map
\[\fa_X : \sPic X(k) \to \Br(k(X)/k)\]
for a smooth projective variety $X$ over $k$. In section~\ref{pairings section} we give a number of different
interpretations of the map $\fa_X$ in terms of pairings and in section~\ref{application section} we use these
to give a number of applications: 
In section~\ref{cyclic pairings} we show that if $X$ is a homogeneous space for
an elliptic curve $E$ of cyclic type (Definition~\ref{cyclic type def}) with respect to the subgroup $T
\subset E$, then we have an exact sequence (generalizing
\cite{Han:RB}), where $E' = E/T$,
\[\xymatrix{
E'(k) \ar[r]^{\phi} & E(k) \ar[r]^-{\fa_X} & \Br(k(X)/k) \ar[r] & 0, 
}\]
and where $\phi$ is dual to the isogeny with kernel $T$.  We further
give an interpretation of the map $\fa_X$ in terms of a natural cup
product.  In section~\ref{perind subsection}, we give a generalization
of a result of Cassels (\cite{Cass:PI}) that the period and index must
coincide for an element of the Tate-Shafarevich group. In
section~\ref{elementary obstruction section} we relate the relative
Brauer group to the elementary obstruction (see
Definition~\ref{elementary obstruction definition}) and prove
\begin{thm*}[\ref{generic_tate}]
Suppose $X$ is a homogeneous space for an elliptic curve $E$ defined over
$k$ and $X(k) = \emptyset$. Then $\Br\big(X_{k(E)}/k(E)\big) \neq 0$ ---
i.e. the relative Brauer group must be nontrivial when one extends scalars
to the function field of $E$.
\end{thm*}
It follows from this (see Corollary~\ref{elementary obstruction}) that if $X(k) \neq \emptyset$ for $X$ as
above, then there exists a field extension $K/k$ such that the elementary obstruction for $X_K$ is nontrivial.
In section~\ref{subsec_app_br} we show that the relative Brauer group of any smooth projective variety $X$ is
always finite when $k$ is local or finitely generated over a prime field (Proposition~\ref{fg finite br}),
however we construct in Theorem~\ref{big_brauer} certain fields $k$, and genus $1$ curves $X/k$, such that the relative
Brauer group $\Br(k(X)/k)$ is infinite.

Finally, in section~\ref{brauer_obstruction} we show that the relative Brauer group may be computed
algorithmically.  The algorithms described in this section have been implemented (in certain cases) as a
Macaulay2 \cite{M2} package ``Relative Brauer,'' freely available at the second author's webpage
\cite{br_cocycle}. This package uses pari as well \cite{Pari}, and may be used to produce examples of relative
Brauer groups for certain homogeneous spaces of elliptic curves defined over $\bQ$. Some of the examples in
this paper were produced using this program.

Although this paper is concerned principally with the computation of relative
Brauer groups, one is often also interested in the more precise question
of how to calculate the index of $\alpha_{k(X)}$ for a general class
$\alpha$ in the Brauer group of $k$. In \cite{LieKra:IR}, the second
author joint with M. Lieblich has shown that this more general problem of
index reduction for a genus $1$ curve may be entirely reduced to the
problem discussed in this paper of computing the relative Brauer group.

\subsection{Definitions and notation}

Throughout the paper we will consider an arbitrary ground field $k$,
and we will denote by $\sep{k}$ a fixed separable closure. We
denote the absolute Galois group $Gal(\sep{k} / k)$ by $G$. Unless
specified otherwise, all cohomology groups should be interpreted as \'etale
cohomology, and in particular, $H^i(k, A)$ coincides with the Group
cohomology $H^i(G, A)$.

\begin{defn}
For a $k$-scheme
$X$, we define the \textit{index} of $X$ to be 
\[\ind(X) = \gcd\{[E:k] | E/k \text{ is a finite field extension and }
X(E) \neq \emptyset\}.\]
\end{defn}
In the case that $X$ is a genus $1$ curve, it follows from
\cite{Lich:PI} that this coincides with the same $\gcd$ taken only over
degrees of separable field extensions.

For a smooth proper variety $X$, we denote its Picard group by
$\Pic{X}$. We define $\sPic{X}$ to be the
sheafification of the $fppf$-presheaf
\[S \mapsto \Pic{X \times_k S}. \] 
This is represented by a $k$-group scheme (see \cite[II.15]{Mur:CF}), which
by abuse of notation, we also denote by $\sPic{X}$. This is projective
in case $X$ is (see \cite[Chap. 9, Thm. 9.5.4]{FGAEx}
We let $\Picnot{X}$ denote the subgroup of $\Pic{X}$ consisting of those
divisor classes which are algebraically equivalent to $0$, and we
recall that the Picard variety of $X$, denoted $\sPicnot{X}$
is the corresponding subscheme of $\sPic{X}$. We recall that
$\sPic{X}$ is an Abelian variety if either $X$ is a curve
\cite[Chap 9, Ex. 9.5.23, pages 289,309]{FGAEx} or an Abelian variety
\cite[III.13]{MumAV}. Since $X$ has a point after
some finite separable field extension \cite[Prop. 10, page
76]{Lang:IAG}, it follows from \cite[Chap. 9, Thm. 9.2.5]{FGAEx} that
$\sPic{X}$ may be considered also as the \'etale sheafification of the
above presheaf.  In particulare may describe the $k$-points of
$\sPic{X}$ as $\sPic{X}(k) = (\Pic{X_{\sep{k}}})^G$.  For every
$G$-fixed element of the Neron-Severi group $\lambda \in \NS(X_{\sep
k})^G$, we may consider the subscheme $\sPicn{\lambda}{X} \subset
\sPic{X}$ of divisors of $X$ in the class $\lambda$. These are
principal homogeneous spaces for the abelian variety $\sPicnot{X}$.
Note that in the case of a curve, we may identify $NS(X) = NS(X_{\sep
k})^G \cong \mbb Z$ by the degree map on divisors, and hence we may
denote classes by integers.  We recall that the collection of all
principal homogeneous spaces for an abelian variety $A/k$ may be
identified with the torsion abelian group $\Het^1(k, A)$, which is
also called the Weil-Chatelet group of $A$.

\begin{defn} \label{period def}
Let $X$ be a smooth projective curve over $k$. We define the
\textit{period} of $X$ to be the order of the class $[\sPicn 1 X]$ in
$\Het^1(k, \sPicnot{X})$.
\end{defn}

For a $k$-scheme $X$, we let $\Br(X)$ denote the Brauer group of
equivalence classes of Azumaya algebras, and for a ring $R$, we write
$\Br(R)$ for $\Br(\Spec(R))$. In the case $X$ is a smooth and
quasi-projective variety, we may identify $\Br(X) = \Het^2(X, \bG_m)$ by the
result of Gabber/deJong (see \cite{deJGab}), and we have an injection
$\Br(X) \hra \Br(k(X))$ (see \cite{Gro:GB1,Gro:GB2,Gro:GB3}). 

\begin{defn} \label{rel_br_def}
Given a morphism of schemes $Y \ra X$, we define the relative
Brauer group, written $\Br(Y/X)$ the kernel of the pullback map $\Br(X)
\ra \Br(Y)$.
\end{defn}
We will frequently abuse notation and write $\Br(Y/R)$ if $X = \Spec(R)$ is
affine, or $\Br(S/R)$ if $Y = \Spec(S)$ and so on.

We recall the definition of the elementary obstruction:
\begin{defn}\cite[Definition~2.2.1]{CTS:DVR2} \label{elementary obstruction definition}
Let $X$ be a smooth geometrically integral variety over $k$, and consider
the exact sequence of Galois modules
\[ 0 \to (\sep k) ^* \to \sep k(X)^* \to \sep k(X)^*/(\sep k) ^* .\]
The elementary obstruction $ob(X)$ is by definition the class of this
extension of $G$ modules in $\operatorname{Ext}_G^1(\sep k(X)^*/(\sep
k)^*, (\sep k)^*)$.
\end{defn}

\begin{Prop}\cite[Proposition 2.2.2(a)]{CTS:DVR2} \label{elementary obstruction vanishes}
Suppose $X$ is a smooth geometrically integral variety over a field $k$ and $X(k) \neq \emptyset$. Then
$ob(X)$ vanishes.
\end{Prop}

\section{The Picard variety, the period and the index} \label{pic period
index}

If $X$ is a smooth projective variety over a field $k$, there is a
well-known natural surjective map
\[ \fa_X : \sPic X (k) \to \Br(k(X)/k) \]
which will be a critical tool in our description of the relative
Brauer group. Particularly important will also be the restriction of
this map to the divisor classes algebraically equivalent to $0$. 

We define the map $\fa_X$ as follows. Consider the short exact
sequence of $G$-modules:
\begin{equation} \label{picard sequence}
\xymatrix{
0 \ar[r] & \sep k (X)^*/(\sep{k})^* \ar[r] & \Div(X_{\sep{k}}) \ar[r] &
\Pic{X_{\sep{k}}} \ar[r] & 0 
}
\end{equation}
Since $\sep k(X)/k(X)$ is a $G$-Galois extension, we may use Hilbert's
Theorem~90 to identify
\[\Hgal^1(G, \sep{k}(X)^*/(\sep{k})^*) \cong \ker\big(\Br(k) \to
\Br(k(X))\big) = \Br(k(X)/k) = \Br(X/k) \]
Using this, the long exact sequence in cohomology from sequence \ref{picard
sequence} gives the desired map:
\[\fa_X : \sPic X \to \Br(k) \]

We now investigate the restriction of this map to the divisor classes
of degree $0$.  Let $\NS(X)$ be the Ner\'on-Severi group
$\Div(X)/\Div^0(X)$, where $\Div^0(X)$ is the group of divisors which
are algebraically equivalent to $0$. The group $\NS(X_{\sep{k}})$
comes with an action of the Galois group $G$, and the short exact
sequence
\[\xymatrix{
0 \ar[r] & \Picnot{X_{\sep{k}}} \ar[r] & \Pic{X_{\sep{k}}} \ar[r] &
\NS(X_{\sep{k}}) \ar[r] & 0 }\]
gives rise to a boundary map
\[\NS(X_{\sep{k}})^G \to \Het^1(k, \sPicnot{X}(\sep k)).\]
Note that $Div(X_{\sep{k}})^G = Div(X)$. We define in $\NS(X_{\sep{k}})^G$ two subgroups $i(X)$ and $p(X)$
as follows:
\begin{gather*}
i(X) = \im\big(\Div(X) \to \NS(X_{\sep{k}})^G\big) = \im\big(\NS(X) \to
\NS(X_{\sep{k}})^G\big) \\
p(X) = \ker\big(\NS(X_{\sep{k}})^G \to \Het^1(k, \sPic{X})\big)
\end{gather*}
These groups have been defined and studied independently by Peter Clark
(see \cite{Clark:PI2}). 

\begin{thm} \label{generalized relative brauer sequence}
Let $X$ be a smooth projective variety over $k$. Then we have exact sequences
\[\xymatrix @C=1.5pc @R=0.5pc {
0 \ar[r] & \Pic X \ar[r] & \sPic X (k) \ar[r]^-{\fa_X} & \Br(X/k) \ar[r] &
0  \\
0 \ar[r] & \Picnot{X} \ar[r] & \sPicnot{X}(k) \ar[r]^-{\fa_X} & \Br(X/k) \ar[r] &
p(X)/i(X) \ar[r] & 0
}\]
\end{thm}

\begin{Remark}
If $X$ is a smooth projective curve this is in essence done by
Lichtenbaum in \cite{Lich}. In this case, we note that $NS(X_{\sep{k}})$ is
simply isomorphic to $\mbb Z$ by associating to every divisor its degree.
In particular, we observe in  this case $i(X) = \ind(X) \mbb Z$. Further,
it follows from \cite{Lich} that the boundary map $\NS(X_{\sep{k}})^G \cong
\mbb Z \to \Het^1(k, \sPic X)$ sends $1$ to the class of the
homogeneous space $\sPicn{1}{X}$, and so we
obtain $p(X) = \per(X) \mbb Z$. We therefore see that the failure of
surjectivity of $\fa_X|_{\sPicnot{X}(k)}$ exactly describes the
obstruction for the period and index of the curve $X$ to coincide.
\end{Remark}

\begin{cor} \label{per_eq_ind}
Suppose $X$ is a curve with $\ind(X) = \per(X)$. Then the map $\fa$ is
surjective.
\end{cor}

\begin{proof}[Proof of Theorem~\ref{generalized relative brauer sequence}]
Identifying $\Pic X$ with the image in $\sPic X (k) = \Pic{X_{\sep k}}^G$
of $\Div(X)$, we obtain from the long exact sequence in cohomology from
sequence \ref{picard sequence}:
\[0 \to \Pic X \to \sPic X (k) \to \Br(k) \to \Br(X) \]
which immediately gives us the first part of the theorem.

For the second part, let $\Br_0(X/k) = \fa_X(\sPicnot X(k))$.
Using the morphism of exact sequences:
\[\xymatrix{
0 \ar[r] & \sep k (X)^*/(\sep{k})^* \ar@{=}[d] \ar[r] & \Div^0(X_{\sep{k}})
\ar[d] \ar[r] & \Picnot{X_{\sep{k}}} \ar[r] \ar[d] & 0  \\
0 \ar[r] & \sep k (X)^*/(\sep{k})^* \ar[r] & \Div(X_{\sep{k}}) \ar[r] &
\Pic{X_{\sep{k}}} \ar[r] & 0 
}\]
We obtain an inclusion of short exact sequences:
\begin{equation}\label{both picards}
\xymatrix @R=1.5pc @C=1.5pc {
0 \ar[r] & \Picnot X \ar@{^{(}->}[d] \ar[r] & \sPicnot X (k)
\ar@{^{(}->}[d] \ar[r] & \Br_0(X/k) \ar@{^{(}->}[d] \ar[r] & 0 \\
0 \ar[r] & \Pic X \ar[r] & \sPic X (k) \ar[r] & \Br(X/k) \ar[r] & 0 \\
}
\end{equation}
Using the long exact sequence associated to the sequence
\[\xymatrix{
0 \ar[r] & \Picnot{X_{\sep{k}}} \ar[r] & \Pic{X_{\sep{k}}} \ar[r] &
NS(X_{\sep{k}}) \ar[r] & 0 \\
}\]
we may identify 
\[p(X) = \coker\big(\sPicnot X \to \sPic X),\] 
and using the sequence
\[\xymatrix{
0 \ar[r] & \Div^0(X_{\sep{k}}) \ar[r] & \Div(X_{\sep{k}}) \ar[r] &
\NS(X_{\sep{k}}) \ar[r] & 0 
}\]
together with the fact that the map $\Div(X) = \Div(X_{\sep k})^G \to
\NS(X_{\sep k})^G$ factors through the surjective map $\Div(X) \to \Pic
X$, we may also identify
\[i(X) = \coker\big(\Picnot X(k) \to \Pic X(k)\big). \]
Therefore, applying the snake lemma to sequence \ref{both picards}, gives
an exact sequence of cokernels:
\[0 \to i(X) \to p(X) \to \frac{\Br(X/k)}{\Br_0(X/k)} \to 0 \]
which gives, via the definition of $\Br_0(X/k)$, 
\[0 \to \Picnot X \to \sPicnot X(k) \to \Br(X/k) \to p(X)/i(X) \to 0\]
as desired.
\end{proof}

\section{Pairings and the map $\fa_X$} \label{pairings section}

In this section we show that the map $\fa_X$ constructed above may be
interpreted in a number of ways in terms of pairings. We will
afterwords use these different interpretations to prove new results
concerning relative Brauer groups in section \ref{application
section}.

Throughout, we will abuse notation and write $\fa_\gamma$ instead of
$\fa_X$ in the case that $X$ is a genus $1$ curve corresponding to a
cohomology class $\gamma \in \Het(k, E)$ for an elliptic curve $E$.

\subsection{The Tate pairing}

\begin{thm}[\cite{Lich:PI}]
Let $A$ be an abelian variety over a field $k$, and consider the Tate
pairing:
\[\left<\ \ ,\ \ \right> : \Het^1(k, A) \times \sPicnot A(k) \to \Br(k).\]
In the case $A$ is an elliptic curve, we have $\left<\gamma, p\right> =
\fa_\gamma(p)$. 
\end{thm}

\begin{proof}
This may be found in \cite[pages~1213-1216]{Lich:PI}. Since it is not
explicitly stated as a theorem in this paper, we note that on page 1213 of
\cite{Lich:PI}, Lichtenbaum defines the pairing due to Tate, on page 1215,
he defines the pairing coming from $\fa$, and in pages 1215-1216 proves
that these coincide. 
\end{proof}

For future reference, let us also recall the definition of the Tate
pairing.  Let $A$ be an abelian variety over $k$.  Denote by $Z(A)$
the group of $0$-dimensional cycles on $A_{\sep{k}}$ of degree $0$,
and by $Y(A)$ the Albanese kernel of $A$ defined by the exact
sequence:
\[0 \to Y(A) \to Z(A) \to A(\sep{k}) \to 0.\]
Let $D \subset A \times \sPicnot A$ be a Poincar\'e divisor and let
$\pi_1, \pi_2$ be the projection maps from $A \times \sPicnot A$ to $A$
and $\sPicnot A$ respectively. We define $Z(A)_D$ (respectively $Y(A)_D$)
to be the subgroup of $Z(A)$ (resp.  $Y(A)$) of elements $\alpha$, such
that $\pi_1^{-1}(|\alpha|)$ transversely intersects $D$, where $|\alpha|$
is the support of $\alpha$.  Note that we have a well defined map
\[\pi_2\left(\pi_1^{-1}(\ \_ \ ) \cap D\right) : Z(A)_D \to \Div(\sPicnot
A_{\sep k}).\]  
Since this agrees with the cycle theoretic map $(\pi_2)_*
\left(\pi_1^*(\alpha \cdot D)\right)$, and since the elements of $Y(A)$
are rationally equivalent to $0$, it follows that the image of this map is
always a principal divisor, giving us actually a map
\[\pi_2\left(\pi_1^{-1}(\ \_\ ) \cap D\right) : Y(A)_D \to \sep
k(A)^*/(\sep k)^*.\] 

To define the Tate pairing, we start with a class $\gamma \in \Het^1(k,
A(\ov k))$, and let $\alpha \in \Het^2(k, Y(A))$ be its image
under the connecting homomorphism. Choose a representative cochain $\ov
\alpha$ for $\alpha$ and choose a Poincar\'e divisor $D$ transversal
to $|\ov \alpha|$. Consider the $\sep k(A)^*/\sep k^*$-valued cocycle
$\beta = \pi_2\left(\pi_1^{-1}(\ov \alpha) \cap D\right)$. Since $A(k)
\neq \emptyset$, the elementary obstruction for $A$ vanishes
(Proposition~\ref{elementary obstruction vanishes}), and consequently,
we may lift the class of $\beta$ to a class $\til \beta \in \Het^2(k,
\sep k(A)^*)$ which turns out to be unramified --- i.e. an element of
$\Het^2(A, \mbb G_m) =
\Br(A)$. By changing $\til \beta$ by a constant class from $\Br(k)$,
we may assume that $\til \beta$ is trivial when specialized to the
identity $0 \in A(k)$. The pairing is then defined by \(\<\gamma, p\> =
\til \beta|_p.\)

\subsection{Pairing via specializations of Brauer classes}

Let $X$ be a smooth projective variety over $k$ with a rational
point $x \in X(k)$. Recall that the group $\Br(X_{\sep k}/X, x)$ is
defined to be the subgroup of $\Br(X_{\sep k}/X)$ consisting of those
Brauer classes $\alpha$ such that the specialization of $\alpha$ at
$x$ is trivial.

\begin{Lem} \label{brauer weil-chatelet isomorphism}
We have an isomorphism
\[\xymatrix @R=.5pc{
H^1(k, \sPic X) \ar[r]^\sim & \Br(X_{\sep k}/X, x)  \\
 \gamma \ar@{|->}[r] & \mc A_\gamma,
}
\]
defined as follows. Consider the exact sequence: 
\begin{equation} \label{pic sequences}
\xymatrix @R=1.5pc @C=1.5pc {
0 \ar[r] & \sep k (X)^*/(\sep k)^* \ar[r] & \Div(X_{\sep{k}}) \ar[r] &
\sPic{X}(\sep k) \ar[r] & 0.
}
\end{equation}
For $\gamma \in \Het^1(k, \sPic{X})$, take the image of $\gamma$ under the
connecting homomorphism in the first exact sequence. This is an
element $\til{\gamma} \in \Het^2(k, \Prin(X)) = \Het^2(k, \sep{k}(X)^*
/ (\sep{k})^*)$. Any lift of this class to $\Het^2(k, \sep k(X)^*)$ will
be an element of the unramified Brauer group of $X$, and we define
$\mc A_\gamma$ to be the unique class
in $\Het^2(k, \sep k (X)^*)$ lifting it which is unramified and is
trivial when specialized to the point $x$.
\end{Lem}

\begin{Rem}
Note that we abuse notation here in that $\mc A_\gamma$ depends on a
particular choice of point $x$. In the case which will be
especially useful to us, $X$ will be an Abelian variety, and in this
case, we will always take the point $x$ to be the identity $0 \in
X(k)$.
\end{Rem}

\begin{Rem}
We will often wish to consider classes $\gamma$ in $H^1(k, \sPicnot
X)$. We abuse notation in this case and write $\mc A_\gamma$ to denote
the Brauer class associated to the image of $\gamma$ in 
$H^1(k, \sPic X)$.
\end{Rem}

We will use this lemma to prove the following alternate formulation of
the Tate pairing:

\begin{thm} \label{alt_tate_pairing}
Let $A$ be an abelian variety over $k$, and $B = \sPicnot A$ its dual.
Let $\alpha \in \Het(k, B(\sep k))$, and choose $X_\alpha$ a homogeneous space
in the class $\alpha$. Then for each $p \in A(k)$, we have:
\[\<\alpha, p\> = \fa_{X_\alpha}(p) = \mc A_{\alpha}|_{p}\]
\end{thm}

We now give a proof of the lemma, followed by a proof of this theorem.

\begin{proof}[Proof of Lemma~\ref{brauer weil-chatelet isomorphism}]
Since $\Div(X_{\sep k})$ has a basis which is permuted by the Galois group,
it follows that $H^1(k, \Div(X_{\sep k})) = 0$ (see for example
\cite[Lemma~12.3]{Sa:LN}). We therefore obtain from (\ref{pic sequences})
an exact sequence
\begin{equation} \label{prin pic}
0 \to \Het^1(k, \sPic{X}) \to \Het^2(k, \sep k(X)^*/(\sep k) ^*) \to
\Het^2(k, \Div(\ov X)).
\end{equation}

Using Proposition~\ref{elementary obstruction vanishes}, we see that since
$X(k) \neq \emptyset$, we have a split exact sequence:
\[0 \to (\sep{k})^* \to \sep{k}(X)^* \to \sep k(X)^*/(\sep k)^* \to 0,\]
and therefore we obtain an exact sequence
\begin{equation} \label{br h2 prin}
0 \to \Br(k) \to \Br(\sep k(X)/k(X)) \to \Het^2(k, \sep k(X)^*/(\sep
k)^*) \to 0.
\end{equation}
The map $\sep{k}(X) \to \Div(X_{\sep{k}})$ induces a map
\[\ram : \Br(\sep k(X)/k(X)) \to \Het^2(k, \Div(X_{\sep{k}}))\]
called the ramification map\footnote{Although this is perhaps not the standard definition of the ramification
map, this definition is shown to be equivalent in \cite{GS:CSAGC} to the standard one, as defined for example
in \cite[Chapter~10]{Sa:LN}. Here our map is given in \cite[Section
6.6]{GS:CSAGC}equation~(6), page 152, and is reformulated on the following
page as being equivalent to (the sum of) the standard ramification
maps over all closed points}, and $\ker(\ram) = \Br(X)$. Let us denote
the kernel of the map
\[\Het^2(k, \sep k(X)^*/(\sep k)^*) \to \Het^2(k, \Div(X_{\sep{k}}))\]
by $\Het^{2, nr}(k, \sep k(X)^*/(\sep k)^*)$.
Using sequence~(\ref{br h2 prin}) we find we have an exact sequence
\[0 \to \Br(k) \to \Br(X_{\sep k}/X) \to \Het^{2, nr}(k, \sep
k(X)^*/(\sep k)^*)
\to 0\]
and an identification $\Het^{2, nr}(k, \sep k(X)^*/(\sep k)^*) =
\Het^1(k, \sPic X)$.
Again using the fact that $X$ has a rational point, we may use
specialization of Brauer classes at the point $x$ to split the map
$\Br(k) \to \Br(X)$, yielding an isomorphism
\[\Br(X_{\sep k}/X, x) \cong \Het^1(k, \sPic X)\]
as desired.
\end{proof}
\begin{Rem}
We note that it follows from this proof that in fact we may express the
relative Brauer group $\Br(X_{\sep k}/X)$ as a product
\[\Br(X_{\sep k}/X) \cong \Br(k) \times H^1(k, \sPic X). \]
Further, in the case that $X$ is a curve, $H^1(k, \sPic X) = H^1(k,
\sPicnot X)$, which follows from the fact that
\[0 \to \sPicnot X(\sep k) \to \sPic X(\sep k) \to \mbb Z \to 0\]
is split exact (since $X(k) \neq \emptyset$), and $H^1(k, \mbb Z) =
0$.
\end{Rem}
  
\begin{proof}[Proof of Theorem~\ref{alt_tate_pairing}]
It is not hard to check that the operation $\pi_B (\blank \cap D)$ gives
a commutative diagram of \textit{sets}:
\begin{equation} \label{tate_diagram}
\xymatrix @R=1pc{
0 \ar[r] & \Prin(B) \ar[r] & \Div^0(B) \ar[r] & A(\sep{k}) \ar[r] & 0 \\
& Y(A)_{D} \ar[u] \ar@{^{(}->}[d] \ar[r] & Z(A)_{D} \ar[u]
\ar@{^{(}->}[d] \ar[r] & A(\sep{k}) \ar@{=}[u] \ar@{=}[d] \\
0 \ar[r] & Y(A) \ar[r] & Z(A) \ar[r] & A(\sep{k}) \ar[r] & 0, } 
\end{equation}
where the top and bottom rows are exact sequences of abelian groups with
Galois action. The pairing of Tate is obtained by the composition of
the connecting homomorphism from
the bottom of diagram (\ref{tate_diagram}) on elements which lie in the
middle row followed by the upwards vertical map, while the pairing using $\mc A$ uses the connecting
homomorphism on the top of the diagram. The result therefore follows
immediately from commutativity of the diagram and equality of the
rightmost terms.
\end{proof}

\subsection{Pairings via torsion points}

For an elliptic curve $E$ over $k$ and an integer $n$ prime to the
characteristic of $k$, we will let $e_n$ denote the Weil
pairing $E[n] \otimes_{\mbb Z} E[n] \to \mu_n$. We let $\delta_n : E(k) \to
\Het^1(k, E[n])$ denote the boundary map from the Kummer sequence:
\begin{equation}\label{kummer}
0 \to E[n] \to E \overset{n}{\to} E \to 0
\end{equation}

\begin{thm} \label{pairings_agree}
Let $E$ be an elliptic curve over $k$, and let 
$p \in E(k)$ and $\gamma \in \Het^1(k, E[n])$, and let $X_\gamma$ be a
homogeneous space in the class $\gamma$. Then
\[\left<\ov \gamma, p\right> = e_n(\gamma \cup \delta_n p) = 
\fa_{X_\gamma}(p) = \mc A_{\gamma}|_p,\]
where $\ov \gamma$ is the image of $\gamma$ in $\Het(k, E)$. 
\end{thm}
\begin{proof}
It follows from \cite[Proposition 9]{Ba}, that $\left<\ov \gamma,
p\right> = e_n(\gamma \cup \delta_n p)$. The remaining assertions
follow from Theorem~\ref{alt_tate_pairing}.
\end{proof}

\section{Applications}\label{application section}

\subsection{Pairings via cyclic isogenies} \label{cyclic pairings}

Let $E$ be an elliptic curve, and suppose we are given a finite Galois
submodule $T \subset E$, we obtain an isogeny, uniquely defined up to
isomorphism \cite[III.\S4, Prop. 4.12]{Sil:EC}:
\begin{equation} \label{isogeny seq}
0 \to T \to E \overset{\phi}{\to} E' \to 0,
\end{equation}
for an elliptic curve $E'$ and a dual isogeny, also unique up to 
isomorphism \cite[III.\S6, Thm. 6.1]{Sil:EC}:
\begin{equation} \label{dual isogeny seq}
0 \to T' \to E' \overset{\phi'}{\to} E \to 0.
\end{equation}

We say that $T \subset E$ cyclic if its points defined over a
separable closure are cyclic as an abstract group. In
other words, $T$ is a closed reduced subscheme of $E$ such that
$T(\sep{k})$ is a cyclic subgroup of $E(\sep{k})$.

\begin{Prop} \label{small_pairing}
Suppose that $T \subset E$ is a cyclic submodule of order $n$ with $n$
prime to the characteristic of $k$, and let $E', T'$ be as above. Then
there is a natural isomorphism of $G$-modules:
\[T \otimes_\bZ T' \cong \mu_n.\]
Further, if $i : T \to E[n]$ is the inclusion and $p = \phi|_{E[n]}$
then this isomorphism is given by $t \otimes t' \mapsto e_n(it,
p^{-1}t')$.
\end{Prop}
\begin{proof}
If $i : T \to E[n]$ is the natural inclusion, it is easy to
see that we obtain a commutative diagram with exact rows:
\begin{equation} \label{isogeny_diagram}
\xymatrix @C=1.5pc @R=1.5pc {
0 \ar[r] & T \ar[r] \ar[d]_i & E \ar@{=}[d] \ar[r]^\phi & E'
\ar[r] \ar[d]^{\phi'} & 0 \\
0 \ar[r] & E[n] \ar[r] \ar[d]_p & E \ar[d]^\phi \ar[r]^n & E
\ar[r] \ar@{=}[d] & 0 \\
0 \ar[r] & T' \ar[r] & E' \ar[r]_{\phi'} & E \ar[r] & 0,}
\end{equation}
which by the snake lemma (applied to the bottom two sequences) gives a
short exact sequence
\begin{equation} \label{finite_duality}
0 \to T \overset{i}{\to} E[n] \overset{p}{\to} T' \to 0.
\end{equation}
where $i$ is induced by the natural inclusion and $p$ is induced by
the map $\phi$.
Consider the Weil pairing:
\[e_n : E[n] \otimes_\bZ E[n] \to \mu_n.\]
Since this pairing is alternating, we have $e_n(T, T) = 0$, and so we
have an induced non-degenerate pairing $T \otimes_\bZ E[n]/T \to
\mu_n$. But by equation (\ref{finite_duality}), we obtain a
non-degenerate pairing $T \otimes_\bZ T' \to \mu_n$. Since this is
clearly an isomorphism ignoring the Galois action, this gives an
isomorphism of Galois modules $T \otimes_\bZ T' \cong \mu_n$ as
desired.
\end{proof}

\begin{defn} \label{cyclic type def}
Let $X$ be a homogeneous space for $E$ of period $n = \per(X)$. We say
that $X$ has cyclic type if its cohomology class in $\Het^1(k, E)$ may be
represented as the image of a cocycle $\gamma \in \Het^1(k, T)$ where $T
\subset E$ is a cyclic submodule of order $n$.
\end{defn}

The following observation shows that cyclic homogeneous spaces are
somewhat special:
\begin{Prop} \label{cyclic perind}
Suppose $X$ is a cyclic homogeneous space for $\E$. Then $\per(X) =
\ind(X)$.
\end{Prop}
\begin{proof}
Choose a particular Galois cocycle $\ov \gamma \in Z^1(G, T)$
representing $X$. Via Galois descent, we may describe $X$ as given by
the curve $\E_{\sep k}$ equipped with the new Galois action $\sigma \cdot p
= \sigma(p) \oplus \ov \gamma(\sigma)$. Let $\phi : \E \to \E'$ be the
isogeny with kernel $T$.
Consider the isogeny $\phi_X : X_{\sep k} \to E'_{\sep k}$ given after the above
identification by $\phi \times_{k} \sep{k}$. We claim that this descends to
give a morphism $X \to \E'$. To see this we need to check
\(\sigma(\phi(p)) = \phi(\sigma \cdot p).\) But since $\ov
\gamma(\sigma) \in T = \ker(\phi)$, we have
\begin{equation*}
\phi(\sigma \cdot p) = \phi(\sigma(p) \oplus \ov \gamma(\sigma)) =
\phi(\sigma(p)) \oplus \phi(\ov \gamma(\sigma)) = \phi(\sigma(p)) =
\sigma(\phi(p)).
\end{equation*}
Now, since we have a $n$ to $1$ \'etale cover $X \to \E'$, the preimage
of the origin in $\E'$ gives a separable point in $X$ of degree $n$ over
$k$, and therefore $\ind(X) | n = \per(X)$. But since $\per(X) |
\ind(X)$ holds for any curve $X$, these must be equal.
\end{proof}

\begin{Lem} \label{cyclic_agrees}
Let $E, E', T, T', \phi, \phi'$ be as above. Let 
\[\delta_{\phi'}: E(k) \to \Het^1(k, T')\]
be the boundary map of the exact
sequence~(\ref{isogeny seq}), and let
\[\delta_n : E(k) \to \Het^1(k, E[n])\]
be the boundary map from the Kummer sequence~(\ref{kummer}). Then we
have for $\gamma \in H^1(k, T)$ and $x \in E(k)$ (and via the identification of
Lemma~\ref{small_pairing}):
\[\gamma \cup \delta_{\phi'} x = e_n(i\gamma \cup \delta_n x )\]
\end{Lem}

Before proving this lemma, we will derive the following consequence, generalizing aspects of the descriptions
of the relative Brauer group obtained in \cite{Han:RB}:
\begin{thm} \label{pairing rel brauer map}
Let $E, E', T, T', \phi, \phi', \delta_\phi'$ be as above. Let 
$X$ be a genus $1$ curve coming from a cocycle $\gamma \in
\Het^1(k, T)$ of order $n$ as above. Then we have a surjective map
\[\xymatrix{
E(k) \ar[r]^-{\fa_X} & \Br(k(X)/k)  
}\]
with $\phi' : E'(k) \to E(k)$ mapping into the kernel of $\fa_X$ and with
$\fa_X$ given by the formula:
$\fa_X(x) = \gamma \cup \delta_{\phi'} x$,
where we identify $T \otimes_\bZ T' \cong \mu_n$ as
in Lemma~\ref{small_pairing}.
\end{thm}
\begin{Rem}
For certain cases of $T$, one may show that the resulting algebras
$\fa_X(x)$ are cyclic as follows:

If $T \cong \mbb Z/n$ then an element $\gamma \in H^1(k, T)$ corresponds
to a cyclic extension of degree $n$ splitting $\fa_X(x)$. Such an
algebra is therefore represented by a cyclic algebra of degree $n$.

In the case that $T$ has odd order and that there exists a quadratic
extension $L/k$ such that $T_L \cong \mbb Z/n$ and such that after
setting $M/L$ to be the cyclic extension corresponding to $\gamma \in
H^1(k_L, T_L)$ we have that $M/k$ is Galois with Galois group
$Gal(M/L) \rtimes Gal(L/k)$. Then in the case that $k$ contains
a primitive $n$'th root of unity, it follows from \cite{RoSa:DA}
that $\fa_X(x)$ may be represented by a cyclic algebra of degree $n$.
More generally, it follows from \cite[Proposition~2.8]{HKRT} that the
same result holds even in the case that $\mu_n$ is contained in any
quadratic extension of $k$.

In the case that $n$ is a prime and $T \cong \mu_n$, $\gamma \in
H^1(k, T)$ corresponds to a Kummer extension splitting $\fa_X(x)$, and
in this case it follows from Albert's result \cite[Theorem 2.11.12,
page 82]{Jac:FDDA} that $\fa_X(x)$ is a cyclic algebra.
\end{Rem}

\begin{proof}[Proof of lemma \ref{cyclic_agrees}]
Choose $x \in \E(k)$ and $\gamma \in \Het^1(k, T)$ and choose $y \in \E(\sep{k})$ with $[n]y = x$. By the commutativity of
diagram (\ref{isogeny_diagram}), $z = \phi(y)$ satisfies $\phi'(z) =
x$. Therefore we have $\delta_{\phi'}(x)(\sigma) = \sigma(z) - z,
\delta_n(x)(\sigma) = \sigma(y) - y$ for $\sigma \in Gal(k)$. 
We now compute:
\begin{align*}
- \gamma \cup \delta_{\phi'}(x)(\sigma, \tau) &=
(\sigma z - z) \otimes \sigma(\gamma(\tau)) \\
&= \left(\sigma(\phi y) - \phi y \right) \otimes \sigma(\gamma(\tau))
\\
&= \phi(\sigma y - y) \otimes \sigma(\gamma(\tau)) \\
&= p\left(\sigma y - y\right) \otimes \sigma(\gamma(\tau)).
\end{align*}
Considering the isomorphism of Proposition~\ref{small_pairing}, this gives:
\begin{align*}
- \gamma \cup \delta_{\phi'}(x)(\sigma, \tau) &=
e_n\left((\sigma y - y) \otimes i\sigma(\gamma(\tau))\right) \\
&= e_n\left((\sigma y - y) \otimes \sigma(i\gamma(\tau))\right) \\
&= - e_n(\delta_n x \cup i\gamma)(\sigma, \tau)
\end{align*}
\end{proof}

\begin{proof}[Proof of Theorem~\ref{pairing rel brauer map}]
We start with the sequence
\[\xymatrix @C=1.5pc{
0 \ar[r] & \Picnot{X} \ar[r] & \sPicnot{X}(k) \ar[r]^-{\fa_X} & \Br(X/k) \ar[r] &
p(X)/i(X) \ar[r] & 0
}\]
of Theorem~\ref{generalized relative brauer sequence}, and identify $E = \sPicnot X$.
Using Theorem~\ref{pairings_agree}, we see that we may represent the map $\fa_X$ as a cup product from the
torsion points of the elliptic curve $E$. Using Lemma~\ref{cyclic_agrees}, we may further interpret this is
coming from the pairing coming from the pair of dual isogenies with kernels $T$ and $T'$. By
Proposition~\ref{cyclic perind}, we have the period and the index of $X$ must coincide, and therefore by
Corollary~\ref{per_eq_ind}, the map $\fa_X$ is surjective.
\end{proof}

\subsection{Period and index} \label{perind subsection}

Since the Hasse principle holds for elements of $\Br(k)$ for $k$ global,
we obtain information about $i(X)$ and $p(X)$ from local data (see section
\ref{pic period index} for the definitions of these invariants). In
particular, we have a simple generalization of a result of (see
\cite{Cass:PI}, \cite{ONeil:PI}, see also \cite{Clark:gen1}), the proof
being an adaptation of the ideas of \cite{ONeil:PI} and \cite{Clark:gen1}
to the present context. This may also be seen as a ``divisorial'' analog
to \cite[Corollary~16]{Olson}.

\begin{thm} \label{indeqexp}
Suppose $k$ is a global field, and let $X/k$ be a smooth projective
variety. Then if $X_v(k_v) \neq \emptyset$ for all but possibly one
valuation $v$ on $k$, then $p(X) = i(X)$.
\end{thm}
\begin{proof}
It suffices to see that these conditions force the relative Brauer group
of $X$ to be trivial. If $\alpha \in \Br(X/k)$, we note that
$\alpha_{k_v} \in \Br(X_{k_v} / k_v)$ for every valuation $v$. By the
hypothesis and by Lemma~\ref{point_no_relbr}, we have $\alpha_{k_v} =
0$. This means that all except possibly one of the Hasse invariants of
$\alpha$ vanishes. But by reciprocity, the sum of the Hasse invariants
is $0$ implying that $\alpha = 0$.
\end{proof}

This immediately implies the result of Cassels \cite{Cass:PI}:
\begin{cor}
If $E$ is an elliptic curve and $[X] \in \Sha(E)$, then $\per(X) =
\ind(X)$.
\end{cor}

\subsection{The elementary obstruction} \label{elementary obstruction section}

The following result says that a nontrivial homogeneous space may always be
detected by its relative Brauer group, at least after extending the
ground field:

\begin{thm} \label{generic_tate}
Suppose $X$ is a homogeneous space for an elliptic curve $E$ defined over
$k$ and $X(k) = \emptyset$. Then $\Br\big(X_{k(E)}/k(E)\big) \neq 0$ ---
i.e. the relative Brauer group must be nontrivial when one extends scalars
to the function field of $E$.
\end{thm}
\begin{proof}
Let $[X]$ denote the class of $X$ in $H^1(k, E)$, and let $\mc A = \mc
A_{[X]}$.
It follows from Theorem~\ref{alt_tate_pairing} that if we consider the generic
point $\eta \in E(k(E))$, then $\fa_X(\eta) = \mc A|_\eta$ is an element
of the relative Brauer group $\Br(X_{k(E)}/k(E))$. Further, since the
isomorphism of Lemma~\ref{brauer weil-chatelet isomorphism} maps the
class $[X]$ to the algebra $\mc A$, it follows that $\mc A$ is a
nontrivial Brauer class.  Since restriction to the generic point gives
an injection $\Br(E) \to \Br(k(E))$, it also follows that $\mc A_\eta$
is nontrivial. Therefore $\Br(X_{k(E)}/k(E)) \neq 0$ as claimed.
\end{proof}

\begin{cor} \label{elementary obstruction}
Let $X$ be a curve of genus $1$ over $k$ and suppose $X(k) \neq
\emptyset$. Then there exists a field extension $K/k$ such that the
elementary obstruction for $X_K$ is nontrivial.
\end{cor}
\begin{proof}
For this, we simply let $K = k(E)$ where $E$ is the Jacobian of $X$. In
this case it follows from Theorem~\ref{generic_tate} that the relative
Brauer group $\Br(K(X)/K)$ is nontrivial. But this implies that the
elementary obstruction must also be nontrivial: arguing by contradiction,
assume $ob(X_K) = 0$. In this case there exists a splitting $\sep K(X)^*
\to (\sep K)^*$, and in particular, the morphism
\[\Br(K) = \Het^2(K, (\sep K)^*) \to \Het^2(K, \sep K(X)^*) = \Br(K(X))\]
is also split injective. But this implies that its kernel, $\Br(K(X)/K)$
must be trivial, yielding a contradiction.
\end{proof}

\subsection{How big is the relative Brauer group?} \label{subsec_app_br}

The following result is a consequence of our results combined with the
theorem of Lang and Neron \cite{LaNe}, and is known by the experts.

\begin{Prop} \label{fg finite br}
Suppose that $X$ is a smooth projective variety defined over a field
$k$ which is local or finitely generated over $\mbb Q$.  Then
$\Br(k(X)/k)$ is finite.
\end{Prop}
\begin{proof}
In the case that $k$ is local this is immediate from
\cite[Theorem~1]{Roq:FF}. If $k$ is finitely generated over its prime
field, we consider the map $\fa_X : \sPicnot{X}(k) \to \Br(k(X) / k)$,
which has finite cokernel by Theorem~\ref{generalized relative brauer
sequence}. Further, since the Brauer group it is enough to show that
$\sPicnot{X}(k)$ is finitely generated. But, since in this case,
$\sPicnot X$ is an Abelian variety over $k$, this follows from the
work of Lang and N\'eron \cite{LaNe}.
\end{proof}

\begin{thm} \label{big_brauer}
Suppose that $E$ is an elliptic curve defined over a field $k$, and
suppose $X$ is a nontrivial homogeneous space over $E$. Then there
exists a field extension $L/k$ such that the relative Brauer group
$\Br(L(X)/L)$ is infinite.
\end{thm}

\begin{Lem} \label{inductive_brauer}
Suppose $X$ is a nontrivial homogeneous space for an elliptic curve
$E$. Then if $L = k(E)$, there is a Brauer class $0 \neq \alpha \in
\Br(X_L/L)$ and an injection $\Br(X/k) \oplus \<\alpha\> \hra
\Br(X_L/L)$.  
\end{Lem}
\begin{proof}
Let $\eta \in E(L)$ be the generic point, and let $\alpha =
\fa_X(\eta)$. By Lemma~\ref{point_splits_proper}, the map $\Br(k) \to
\Br(L)$ is injective. By Theorem~\ref{alt_tate_pairing} and
Theorem~\ref{generic_tate}, the element $\alpha$ is the restriction
of the Brauer class $\mc A(X) \in \Br(E, 0_{E}) \subset \Br(E)$ to the
generic point $\eta \in E$, and $\alpha \neq 0$. 

Using Lemma~\ref{brauer weil-chatelet isomorphism}, we may write
$\Br(E, 0_{E}) \oplus \Br(k) = \Br(E) \hra \Br(L)$. In particular,
since $\alpha \in \Br(E, 0_{E})$ and $\Br(X/k) \subset \Br(k)$, the
groups $\<\alpha\>$ and $\Br(X/k)$ do not intersect considered as
subgroups of $\Br(X_L/L) \subset \Br(L)$. In particular, we obtain an
injection $\Br(X/k) \oplus \<\alpha\> \hra \Br(X_L/L)$ as claimed.
\end{proof}

\begin{proof}[Proof of Theorem~\ref{big_brauer}]
Suppose $X$ is a homogeneous space of index $n$. We will begin by
reducing to the case that $n$ is prime. This is not essential for the
result, but helps the exposition of the proof. Let $p$ be a prime
divisor of $n$, and let $F'/k$ be a prime to $p$ closure of $k$. Note
that $X(F') = \emptyset$ still with $ind(X_F') = p^k$. Since every field
extension has degree a power of $p$, we have 
\[\ind(X_{F'}) = \gcd\{[E:F'] | X(E) \neq \emptyset\} = \min\{[E:F' |
X(E) \neq \emptyset\}.\]
Consequently there is a field $E/F'$ of degree $p^k$ with $X(E) \neq
\emptyset$. Since $F'$ is prime to $p$-closed, it follows that there is
an intermediate field extension $F' \subset F \subset E$ with $[F : F']
= p^{k-1}$. Consequently, $X_F$ has index exactly $p$.

We will construct a chain of field extensions of $F$, 
\[F = L_0 \subset L_1 \subset L_2 \subset \cdots, \]
such that $(\bZ/p)^i \subset \Br(X_{L_i}/L_i)$, and
$\Br(X_{L_{i-1}}/L_{i-1})$ injects into $\Br(X_{L_i} / L_i)$. Assuming
that this has been done for $i-1$, we define $L_i$ to be the function
field $L_{i-1}(X)$. By Lemma~\ref{inductive_brauer}, there is an $\alpha
\in \Br(L_i(X) / L_i)$ such that $\Br(L_{i-1}(X)/L_{i-1}) \oplus
\<\alpha\> \hra \Br(L_i(X) / L_i)$. Since the index of $X_{L_i}$ is $p$,
$\<\alpha\> \cong \bZ/p$, and the induction step follows from the fact
that $(\bZ/p)^{i-1} \subset \Br(L_{i-1}(X)/L_{i-1})$.

Let $L = \underset\to\lim L_i = \cup_i L_i$. Clearly the natural
restriction map $\Br(L_i) \to \Br(L)$ maps $\Br(L_i(X)/L_i)$ to
$\Br(L(X)/L)$. I claim that this map is injective. Arguing by
contradiction, let us suppose there is an $\alpha \in \Br(L_i(X)/L_i)$
with $\alpha_L = 0$. If $A$ is a central simple algebra in the class of
$\alpha$, then this says that the Severi-Brauer variety $SB_A$ has an
$L$-rational point \cite{Sa:LN}. With respect to some projective
embedding of the variety $SB_A$, this point has a finite number of
coordinates, which must therefore lie in some field $L_j$, for a
sufficiently large integer $j$. But this means $SB_A(L_j) \neq
\emptyset$ and so $L_j$ splits $A$. This implies that $\alpha_{L_i} =
0$. But this contradicts the injectivity of $\Br(X_{L_{i-1}}/L_{i-1})
\to \Br(X_{L_i}/L_i)$.

We therefore have $\cup \Br(L_i(X) / L_i) \subset \Br(L(X)/L)$, which
implies $(\bZ/p)^\infty \subset \Br(L(X)/L)$ as desired.
\end{proof}

\section{An explicit description of $\fa_X$ for genus $1$ curves} \label{brauer_obstruction}

Let $\E$ be an elliptic curve over $k$ given on an affine patch by the
equation:
\begin{equation*}
y^2 + a_1 xy + a_3 y = x^3 + a_2 x^2 + a_4 x + a_6
\end{equation*}

Let $L/k$ a $G$-Galois extension (which is no longer assumed to be the
entire absolute Galois group), and let $\gamma \in Z^1(G, \E(L))$ be a
1-cocycle (crossed homomorphism) representing a homogeneous space
$X/k$ for $\E/k$. That is to say, as $G$ varieties, $X_L$ is
isomorphic to $\E_L$ with the Galois action $\sigma^\gamma =
\oplus_{\gamma(\sigma)} \circ \sigma$, where by $\oplus_p$ we mean the
automorphism of the elliptic curve given by addition by $p \in
\E$. This means, for example, that for $p \in \E(L)$, we have
$\sigma^\gamma(p) = \gamma(\sigma) \oplus \sigma(p)$. With this in
mind, we represent points in $X(L)$ by points in $\E(L)$, just with a
different $G$-module structure.

For a function $f \in L(X)$, where $X$ is a $G$-variety, we have an
action of $\sigma \in G$ on $f$ by $\sigma(f) = \sigma \circ f \circ
\sigma^{-1}$ where the $\sigma^{-1}$ is the action on $X$ and the
$\sigma$ is induced by the action on $L$. In particular, if we
identify $L(X)$ with $L(\E)$ with a twisted action, we may write our
action of $\sigma \in G$ on $f \in L(X)$ via
\begin{equation} \label{function_action}
\sigma^\gamma(f)(p) = \sigma \circ f \circ (\sigma^{\gamma})^{-1} (p) =
\sigma \circ f \circ \sigma^{-1}(p \ominus \gamma(\sigma))
\end{equation}

\subsection{Computations} \label{computations}

The goal of this section is to explicitly describe the map $\fa_X :
\E(k) \to \Br(X/k)$ described above in Theorem~\ref{generalized
relative brauer sequence}. Given an element
$p \in \E(k)$, this works in the following steps:
\begin{enumerate}
\item \label{pic_step} Represent $p$ as an element in $(\Picnot{X_{L}})^G$.
\item \label{div_step} Pull this element back to a element in $D_p \in
\Div^0(X_{L})$.
\item \label{cobound_step} Compute the coboundary $\del D_p$ as a
1-cocycle with values in $\Div^0(X_{L})$.
\item \label{prin_step} Realize these values as lying in principal
divisors on $X$ - i.e. for each $p, \sigma$, find a function $f_{p,
\sigma} \in L(X)$ whose divisor is $\del D_p(\sigma)$. This gives a
1-cochain $f_p(\sigma) = f_{p, \sigma}$.
\item \label{function_step} Let $\widetilde{c}_p = \del f_p$, and note
that we may consider this as a 2-cocycle with values in $L^*$ (i.e.
values are constant). That is to say, choosing $q \in X(L)$, we have a
2-cocycle $c_p(\sigma, \tau) = \widetilde{c}_p(\sigma, \tau)(q)$. This
is our Brauer group element.
\end{enumerate}

It will be useful to have an explicit way to show that certain
divisors are principal. We begin with the following definition:
\begin{defn} \label{linears}
Suppose $p, q \in \E(L)$. Define the function $l_{p, q} \in L(\E)$
in the following way:
\begin{itemize}
\item
if $p = q = \infty$, then $l_{p,q} = 1$.
\item
if $p \neq q = \infty$, $p = (x_1, y_1)$ then $l_{p,q} = x - x_1$.
\item
if $p, q \neq \infty$, $p = q = (x_1, y_1)$, then 
\begin{align*}
l_{p,q} = &(y - y_1)(2 y_1 + a_1 x_1 + a_3) \\ 
&- (x - x_1)(3 x_1^2 + 2 a_2 x_1 + a_4 - a_1 y_1).
\end{align*}
\item
if $p, q \neq \infty$, $p \neq q$, $p = (x_1, y_1), q = (x_2, y_2)$,
then 
$$l_{p,q} = (y_2 - y_1)x - (x_2 - x_1)y + x_2 y_1 - x_1 y_2.$$
\end{itemize}
\end{defn}

\begin{Lem} \label{explicit_sum}
Let $p_1, p_2 \in \E(L)$, and let $q = p_1 \oplus p_2$. Then
$$\left(\frac{l_{p_1, p_2}}{l_{q, \ominus q}}\right) = p_1 + p_2 - q -
0_\E.$$
\end{Lem}
\begin{proof}
This is a routine verification. Note that $l_{p,q}$ is the equation of
a line in $\bA^2$ which passes through the points $p, q$.
\end{proof}

We now go through the above steps in sequence:

\begin{enumerate}
\item[\ref{pic_step}.] 

For $p \in \E(k)$, we represent $p$ by the class of the divisor $p -
0_\E \in (\Picnot{X_L})^G$

\item[\ref{div_step}.] 

This is done already. The divisor is $D_p = p - 0_\E$.

\item[\ref{cobound_step}.] 

Set $d_p = \del(D_p)$. explicitly, we have:
\begin{align*}
d_p(\sigma) &= \sigma^\gamma(p - 0_\E) - (p - 0_\E) \\ 
&= \gamma(\sigma) \oplus \sigma(p) - (\gamma(\sigma) + p) + 0_\E.
\end{align*} 
Since $p \in \E(k)$, $\sigma(p) = p$, and so we have:
$$d_p(\sigma) = \gamma(\sigma) \oplus p + 0_\E - \gamma(\sigma) - p.$$

\item[\ref{prin_step}.]
By Lemma~\ref{explicit_sum}, if we set 
$$f_{p, \sigma} = \frac{l_{\gamma(\sigma) \oplus p, \ominus
\gamma(\sigma) \ominus p}}{l_{\gamma(\sigma), p}},$$ we have 
\begin{equation} \label{br_func}
(f_{p, \sigma}) = \gamma(\sigma) \oplus p + 0_\E - \gamma(\sigma) - p =
d_p(\sigma),
\end{equation}
then $f_p(\sigma) = f_{p, \sigma}$ then gives a 1-cochain with values in
$L(X)$.

\item[\ref{function_step}.] 
Let $\widetilde{c}_p = \del f_p$. 
By standard arguments, we in fact know that this function has values in
$L$.  Explicitly we have:
\[c_p(\sigma, \tau) = \frac{f_{p, \sigma} \sigma^\gamma f_{p,
\tau}}{f_{p, {\sigma \tau}}}  = \frac{(f_{p, \sigma})
(\ominus_{\gamma(\sigma)} f_{p, \tau}^{\sigma})}{f_{p, {\sigma \tau}}}\]
\end{enumerate}

\begin{Rem}
The only relevant issue about the functions $f_{p,\sigma}$ is that their
associated principal divisor is described as in equation \ref{br_func}.
In particular, we may change the functions $f_{p,\sigma}$ by any constants
and get an equivalent 2-cocycle describing $c_p$.
\end{Rem}

\begin{thm} \label{presentation}
Suppose $p \in \E(k)$. Then with the above notation, 
\begin{gather*}
\fa_p = (L/k, G, c_p), \\
\text{where $c_p$ is given as:   }
c_p(\sigma, \tau) = \frac{(f_{p, \sigma})(\ominus_{\gamma(\sigma)}f_{p,
\tau})}{f_{p, {\sigma \tau}}}
\end{gather*}
\end{thm}

We gather these facts in the following application:

\begin{cor} \label{computation_cor}
Suppose $X$ is a genus $1$ curve with Jacobian $\E$, and suppose
$\ind(X) = \per(X)$. Let $L/k$ be a Galois extension with group $G$, and
suppose that $X(L) \neq \emptyset$. Then the relative Brauer group
$\Br(X/k)$ is given by
\begin{gather} \label{cocycle_expression}
\{[(L/k, G, c_p)] | p \in \E(k)\} \text{, where} \\ c_p(\sigma, \tau) =
\frac{(f_{p, \sigma}) (\ominus_{\gamma(\sigma)} f_{p, \tau})}{f_{p,
{\sigma \tau}}},
\end{gather}
and where the function $f_{p, \sigma} \in L(\E)$ is defined by the
expression:
$$f_{p, \sigma} = \frac{l_{\gamma(\sigma) \oplus p, \ominus
\gamma(\sigma) \ominus p}}{l_{\gamma(\sigma), p}},$$ and the functions
$l_{p_1, p_2}$ are given in definition \ref{linears}.
\end{cor}

\begin{Rem} \label{rational_cocycle}
In the case that the cocycle $\gamma \in \Het^1(L, \E(L))$ has values in
$\E(k)$, we may simplify our expression for the cocycle $c_p$ above by
noting that the functions $f_{p, \sigma}$ are in $k(\E)$ (using their
explicit description above), and hence are Galois invariant. We may
therefore write in this case:
\[c_p(\sigma, \tau) = \frac{f_{p, \sigma} (\ominus_{\gamma(\sigma)} f_{p,
\tau})}{f_{p, {\sigma \tau}}}\]
\end{Rem}

\begin{Rem} \label{cyclic_relbrauer}
In the case that $G = \<\sigma \mid \sigma^m\>$ is a cyclic group with
generator $\sigma$ of order $m$, and the values of the cocycle are in
$\E(k)$, we may simplify our formula for $c_p$ significantly.  Using the
fact that $c_p$ is a cocycle with values in $k^*$, One may check explicitly
that in this case $c_p$ is cohomologous to $c_p'$ where for $0 \leq i,j <
n$ we have:
\[c_p'(\sigma^i, \sigma^j) = \left\{
\begin{matrix}
1 && i + j < m \\
c_p(1,1) c_p(2,1) \cdots c_p(m-1,1) && i + j \geq m 
\end{matrix}
\right. \]
In particular, the central simple algebra represented by $c_p$ is the
cyclic algebra \[(L, \sigma, c_p(1,1) c_p(2,1) \cdots c_p(m-1,1)).\]
\end{Rem}

\subsection{Examples}

\subsubsection{A hyperelliptic curve}
We demonstrate this formula for the relative Brauer group by
reproducing an example of I. Han \cite{Han:RB}. Let $X$ be the
hyperelliptic genus $1$ curve given by the affine equation $y^2 = ax^4
+ b$, and suppose $X(k) = \emptyset$. The Jacobian $\E$ of this curve
is given by the Wierstrauss equation $y^2 = x^3 - 4abx$. Let $L =
k(\beta)$, where $\beta^2 = b$. We have $(0, \pm \beta) \in X(L)$, and
so the index of $X$ is 2 and hence is also equal to the period of
$X$. Let the order $2$ group $G = <\sigma \mid \sigma^2>$ be the Galois
group of $L/k$.

Note that the (non-identity) 2-torsion points of $\E$ are exactly
those points with $y$ coordinate $0$. These are:
\begin{gather*}
t_0 = (0, 0), \ \  
t_+ = (0, 2\lambda), \ \
t_- = (0, -2\lambda),
\end{gather*}
where $\lambda^2 = 4ab$. In particular, we have three possibilities:
$\lambda \in k$, $\lambda \in L \setminus k$, or $\lambda \not \in L$. We
will assume the second possibility holds : i.e. $ab \in (L^*)^2
\setminus (k^*)^2$.

\begin{Lem}
Define a 1-cocycle $\gamma \in \Hgal^1(G, \E(L))$ via $\gamma(id) = 1,
\gamma(\sigma) = t_0$. Then $\gamma$ corresponds to $X$ viewed as a
homogeneous space over $\E$.
\end{Lem}
\begin{proof}
This is exactly \cite[Example 3.7, pages 293-295]{Sil:EC}.
\end{proof}

We will consider the case where $rk(\E) = 0$. In this case since the
image of $\fa$ is entirely 2-torsion, $im(\fa) = \fa(\E(k)_2)$, where
$\E(k)_2$ is the 2-power torsion part of $\E(k)$. On the other hand,
since the only $k$-rational non-identity 2-torsion point is $t_0$, if
there are other points in $(k)_2$, there must be at least a 4-torsion
point. An explicit computation quickly shows however, that this would
contradict the fact that $ab$ is not a square (the line from $t_0$ to
the 4-torsion point would have to be tangent at the 4-torsion point,
and therefore its slope would have to be a fourth root of $16ab$).
Consequently, we have \[\Br(X/k) = \<\fa_{t_0}\>.\]

By Corollary~\ref{computation_cor}, we have
\[A = \fa_{t_0} = (L/k, G, c_{t_0}),\]
with $c = c_{t_0}$ described as above. For this example we compute explicitly:
\[ f_{t_0, id} = \frac{l_{0_\E \oplus t_0, \ominus \gamma(id) \ominus
t_0}}{l_{\gamma(\sigma), t_0}} = \frac{l_{t_0, t_0}}{l_{0_\E, t_0}} =
\frac{x}{x} = 1, \]
and
\[ f_{t_0, \sigma} = \frac{l_{\gamma(\sigma) \oplus t_0, \ominus
\gamma(\sigma) \ominus t_0}}{l_{\gamma(\sigma), t_0}} = \frac{l_{0_\E,
0_\E}}{l_{t_0, t_0}} = \frac{1}{x}. \]
Using $q = t_+$ in the formula \ref{cocycle_expression} from
Corollary~\ref{computation_cor}, and using the fact that $\sigma(t_+)
= t_-$, we have:
\begin{gather*}
c(id, id) = c(id, \sigma) = c(\sigma, id) = 1, \\
c(\sigma, \sigma) = -\frac{1}{4ab}
\end{gather*}
Therefore, the relative Brauer group $\Br(X/k)$ is generated by the algebra
\begin{equation*}
A = (L/k, \sigma, -\frac{1}{4ab}) = (L/k, \sigma, -ab) = (b, -ab)_{-1}
 = (b, a)_{-1} + (b, -b)_{-1} = (a, b)_{-1}
\end{equation*}
and so $\Br(X/k) = \{1, (a, b)_{-1}\}$.

\subsubsection{A curve of index $5$}

Let $\E$ be the elliptic curve over $\bQ$ given by the equation:
$$y^2 + y = x^3 - x^2 - 10x - 20.$$ One may check (for example, using
Pari \cite{Pari}) that the torsion subgroup of this curve is $\bZ /
5\bZ$, generated by $g = (5, 5)$, and that its rank is $0$. Let $L/\bQ$
be any cyclic degree $5$ extension, say $G = Gal(L/\bQ) = \<\sigma \mid
\sigma^5\>$. Let $\gamma \in \Hgal^1(G, \E(L))$ be given by
$\gamma(\sigma) = g$. The element $\gamma$ corresponds to a
homogeneous space $X$ for $\E$, which is isomorphic to $\E$ if and
only if $\gamma$ is trivial. We will show:

\begin{Ex}
The relative Brauer group $\Br(X/\bQ)$ of the curve $X$ is cyclic
generated by the cyclic algebra $(L/\bQ, \sigma, 11)$. In particular,
if $11$ is not a norm from $L$, then $X$ is not split.
\end{Ex}
\begin{proof}
This follows from direct computation with the above formulas. In
particular, since for our curve, $\E(\bQ) = \<g\>$, one need only check
the 2-cocycle in the image of the point $g$.

One may check that for our curve, we may use the functions $f_{g, \sigma^i}$
given by:
\[\xymatrix @C=-.9pc @R=0.0pc {
f_{g, id} = 1 & & f_{g, \sigma} = \frac{x - 16}{5x-y-20} & &
f_{g, \sigma^2} = \frac{x - 16}{6x+y-35} \\ 
& f_{g, \sigma^3} = \frac{x - 5}{-5x+y+20} & & 
f_{g, \sigma^4} = \frac{1}{x - 5}
}\]
Using remark \ref{rational_cocycle}, we may
express our cocycle $c_g$ as
$$c_g(\sigma^i, \sigma^j) = \frac{f_{g, \sigma^i} (\ominus_{g^i} f_{q,
\sigma^j})}{f_{q, \sigma^{i+j}}}.$$
With the aid of computational software (\cite{Pari}, \cite{M2}), we may
determine the functions $\ominus_{g^i} f_{q, \sigma^j}$. In particular,
we have:
\[\xymatrix @C=0.0pc @R=0.5pc{
\ominus_{g} f_{g, \sigma} = \frac{5x - y - 20}{6x + y - 35} &
\ominus_{g^2} f_{g, \sigma} = \frac{5x + y - 19}{6x - y - 36} \\
\ominus_{g^3} f_{g, \sigma} = \frac{x - 5}{-5x - y + 19} &
\ominus_{g^4} f_{g, \sigma} = \frac{5 - x}{11}
}\]
and by using \ref{cyclic_relbrauer},
we may find a cohomologous cocycle
$$c_g'(\sigma^i, \sigma^j) = \left\{
\begin{matrix}
1 && \text{ if } i + j < 5 \\
1/11 && \text{ if } i + j \geq 5
\end{matrix}
\right.$$ 
\end{proof}

\subsubsection{An example with noncyclic relative Brauer group}

Using the computer package \cite{br_cocycle}, we may construct other
interesting examples by using curves of rank $0$ with interesting
torsion subgroups. The following is the result of output from
this program: 

Let $\E$ be the elliptic curve defined over $\bQ$ by the affine equation
\[y^2 + xy + y = x^3 + x^2 -10x -10\]
Let $L/\bQ$ be a cyclic Galois extension with generator $\sigma$ of
order $4$. Given a torsion point $t \in \E(k)$ of order $n$ dividing
$4$, we may use it to define a homomorphism $Gal(L/k) \to \E(k)$ by
sending the generator $\sigma \in Gal(L/\bQ)$ to the torsion point $t$.
Via the map $Hom(Gal(L/\bQ),\E(k)) \to \Het^1(\bQ, \E)$, this defines a
principal homogeneous space $X_t$.

For the elliptic curve $\E$, pari/gp tells us that $\E$ is rank $0$ with
torsion subgroup generated by the points $(8,18)$ of order $4$ and
$(-1,0)$ of order $2$.

For the homogeneous space defined by $t = \ominus(8,18)$, 
the relative Brauer group is isomorphic to $\bZ/4 \times \bZ/2$,
generated by the cyclic algebras: 
\begin{align*}
\fa_{X_t}(8,18) &= (L/k, \sigma, 405) = (L/k, \sigma, 5), \\
\fa_{X_t}(-1,0) &= (L/k, \sigma,-81) = (L/k, \sigma, -1).
\end{align*}

\appendix

\section{Rational points over function fields}

The following lemma is due to Nishimura \cite{Nish}:
\begin{Lem} \label{point_splits_proper}
Suppose $X$ and $Y$ are schemes over $k$ such that $Y$ is proper, $Y(k)
= \emptyset$, $X$ is irreducible and $x \in X(k)$ is a smooth point.
Then $Y(k(X)) = \emptyset$.
\end{Lem}
\begin{proof}
The proof proceeds by induction on $\dim(X)$. Suppose that $Y(k(X)) \neq
\emptyset$. Then there is a rational morphism $\phi : X \dra Y$. If
$\dim(X) = 0$ then $x = X$ and $\phi$ gives an element of $Y(k)$,
contradicting our hypothesis $Y(k) = \emptyset$. For the general
induction step, let $\til{X}$ by the blowup of $X$ at the point $x$ and
let $E \subset \til{X}$ be the exceptional divisor. Since 
the map $\phi$ may be defined in a set of codimension at
least $2$ and in particular, by restricting this morphism to $E$, we
obtain a rational map $E \dra Y$. Since $E \cong \bP^{\dim(X) - 1}$, $E$
contains a smooth $k$-point and is irreducible.  Therefore by letting
$E$ take the role of $X$, the induction hypothesis implies $Y(k) \neq
\emptyset$, contradicting our hypothesis and completing the proof.
\end{proof}

One consequence of this fact is that index of projective varieties is
not changed by such field extensions:
\begin{cor} \label{point_preserves_index}
Suppose $Y$ is a projective variety, and $X$ is a variety with a smooth
rational $k$-point. Then $\ind Y = \ind Y_{k(X)}$.
\end{cor}
\begin{proof}
For a positive integer $n$, let $Y^{[n]}$ be the Hilbert scheme of
$n$ points on $Y$. Since this is a projective scheme,
Lemma~\ref{point_splits_proper} tells us that $Y^{[n]}(k) \neq
\emptyset$ if and only if $Y^{[n]}(k(X)) \neq \emptyset$. Since the
index of $Y$ may be thought of as the $\gcd$ of the set of positive
integers $n$ such that $Y^{[n]}$ has a rational point, we obtain $\ind
Y = \ind Y_{k(X)}$.
\end{proof}

In particular, we obtain:

\begin{cor} \label{point_no_relbr}
Suppose $X$ is a scheme defined over $k$ with a smooth point $x \in
X(k)$.  Then the restriction map $\Br(k) \to \Br(k(X))$ is injective.
\end{cor}
\begin{proof}
Recall that for a central simple algebra $A$ over a field $F$, $A$ is
split if and only if the associated Severi-Brauer variety $SB_A$ has an
$F$-point \cite{Sa:LN}.  Therefore, if $A$ is a central simple $k$
algebra with $[A] \neq 0$ in $\Br(k)$, then the variety $SB_A$ has no
$k$-points. If the algebra is split by $k(X)$ then this implies $SB_A$
does have a point over $k(X)$. Since $SB_A$ is a proper variety, this
would contradict Lemma~\ref{point_splits_proper}.  Therefore we must
have $SB_A(k(X)) = \emptyset$ and $[A]_{k(X)} \neq 0$.
\end{proof}

\nocite{Haile:CABC}
\nocite{Haile:CABCS}
\nocite{Haile:CACSFF}
\nocite{HaileTess}
\nocite{HodTess}

\bibliographystyle{alpha}
\bibliography{citations}

\end{document}